\input amstex
\input epsf
\documentstyle{amsppt}
\magnification \magstep1
\hsize=6.1truein
\vsize=8.5truein
\hfuzz 7truept

\def\ls{\leqslant}
\def\gs{\geqslant}
\def\flm{f_{l,m}}
\def\plm{P_{l,m}}
\def\flma{f_{l,m}^a}
\def\Flma{F_{l,m}^a}
\TagsOnRight

\topmatter
\title
Avoiding maximal parabolic subgroups of $S_k$
\endtitle
\author Toufik Mansour$^*$ and Alek Vainshtein$^\dag$ \endauthor
\affil $^*$ Department of Mathematics\\
$^\dag$ Department of Mathematics and Department of Computer Science\\ 
University of Haifa, Haifa, Israel 31905\\ 
{\tt tmansur\@study.haifa.ac.il},
{\tt alek\@mathcs.haifa.ac.il}
\endaffil

\abstract
We find an explicit expression for the generating function of the number 
of permutations in $S_n$ avoiding a subgroup of $S_k$ generated by all
but one simple transpositions. The generating function turns out to be
rational, and its denominator is a rook polynomial for a rectangular
board.
\medskip
\noindent {\smc 2000 Mathematics Subject Classification}: 
Primary 05A05, 05A15; Secondary 05E35, 33C45
\endabstract

\leftheadtext{Toufik Mansour and Alek Vainshtein}
\endtopmatter


\document
\heading 1. Introduction and Main Result\endheading

Let $[p]=\{1,\dots,p\}$ denote a totally ordered alphabet on $p$ letters, and 
let $\alpha=(\alpha_1,\dots,\alpha_m)\in [p_1]^m$, 
$\beta=(\beta_1,\dots,\beta_m)\in [p_2]^m$. We say that $\alpha$ is
{\it order-isomorphic\/} to $\beta$ if for all $1\ls i<j\ls m$ one has
$\alpha_i<\alpha_j$ if and only if $\beta_i<\beta_j$. For two permutations
$\pi\in S_n$ and $\tau\in S_k$, an {\it occurrence\/} of $\tau$ in $\pi$
is a subsequence $1\ls i_1<i_2<\dots<i_k\ls n$ such that $(\pi_{i_1},
\dots,\pi_{i_k})$ is order-isomorphic to $\tau$; in such a context $\tau$ is
usually called the {\it pattern\/}. We say that $\pi$ {\it avoids\/} $\tau$,
or is $\tau$-{\it avoiding\/}, if there is no occurrence of $\tau$ in $\pi$. 
Pattern avoidance proved
to be a useful language in a variety of seemingly unrelated problems, from
stack sorting \cite{Kn, Ch.~2.2.1} to singularities of Schubert varieties 
\cite{LS}.
A natural generalization of single pattern avoidance is {\it subset 
avoidance\/};
that is, we say that $\pi\in S_n$ avoids a subset $T\subset S_k$ if
$\pi$ avoids any $\tau\in T$. 
A complete study of subset avoidance for the case $k=3$ is
carried out in \cite{SS}. For $k>3$ situation becomes more complicated,
as the number of possible cases grows rapidly. Recently, several authors
have considered the case of general $k$ when $T$ has some nice algebraic 
properties. Paper \cite{BDPP} treats the case when $T$ is the 
centralizer of $k-1$ and $k$ under the natural action of $S_k$ on
$[k]$ (see also Sec.~3 for more detail). In 
\cite{AR}, $T$ is a Kazhdan--Lusztig cell of $S_k$, or, equivalently, the
Knuth equivalence class (see \cite{St, vol.~2, Ch.~A1}). 
In this paper we consider 
the case when $T$ is a maximal parabolic subgroup of $S_k$.

Let $s_i$ denote the simple transposition interchanging $i$ and $i+1$. Recall 
that a subgroup of $S_k$ is called {\it parabolic\/} if it is generated 
by $s_{i_1},\dots, s_{i_r}$. A parabolic subgroup of $S_k$ is called
{\it maximal\/} if the number of its generators equals $k-2$. We denote
by $P_{l,m}$ the (maximal) parabolic subgroup of $S_{l+m}$ generated
by $s_1,\dots,s_{l-1},s_{l+1},\dots,s_{l+m-1}$, and by $f_{l,m}(n)$ 
the number of permutations in $S_n$ avoiding all the patterns in $\plm$.
In this note we find an explicit expression for the generating function
of the sequence $\{\flm(n)\}$.

To be more precise, we prove the following more general result. Let
us denote
$\sigma=s_1s_2\dots s_{k-1}$, that is, $\sigma=(2,3,\dots,k,1)$,
and let $a$ be an integer, $ 0\ls a\ls k-1$ (here and in what follows $k=l+m$).
We denote by $\flma(n)$ the number of permutations in $S_n$ avoiding the
left coset $\sigma^a\plm$; in particular, $f_{l,m}^0(n)$ coincides with
$\flm(n)$. Let $\Flma(x)$ denote the generating function of $\{\flma(n)\}$,
$$
\Flma(x)=\sum_{n\gs0}\flma(n)x^n.
$$
Recall that the {\it Laguerre polynomial\/} $L_n^\alpha(x)$ is given by
$$
L_n^\alpha(x)=\frac 1{n!}e^xx^{-\alpha}\frac{d^n}{dx^n}
\left(e^{-x}x^{n+\alpha}\right),
$$
and the {\it rook polynomial\/} of the rectangular $s\times t$ board
is given by
$$
R_{s,t}(x)=s!x^sL_s^{t-s}(-x^{-1})
$$
for $s\ls t$ and by $R_{s,t}(x)=R_{t,s}(x)$ otherwise (see \cite{Ri, Ch.~7.4}).

\proclaim{Main Theorem} Let $\lambda=\min\{l,m\}$, $\mu=\max\{l,m\}$, then
$$
\Flma(x)R_{l,m}(-x)=\sum_{r=0}^{\lambda-1}x^rr!
\sum_{j=0}^r(-1)^j\frac{\binom lj\binom mj}{\binom rj}+
(-1)^\lambda x^\lambda\lambda!\sum_{r=0}^{\mu-\lambda-1}x^rr!
\binom {\mu-r-1}\lambda,
$$
or, equivalently,
$$
\Flma(x)=\sum_{r=0}^{k-1}x^rr!-\frac{(-1)^\lambda x^\mu}
{\lambda!L_\lambda^{\mu-\lambda}(x^{-1})}\sum_{r=0}^{\lambda-1}
(k+r)!x^r\sum_{j=r+1}^\lambda(-1)^j\frac{\binom lj\binom mj}{\binom {k+r}j},
$$
where $k=l+m=\lambda+\mu$.
\endproclaim

The proof of the Main Theorem is presented in the next section.

As a corollary we immediately get the following result (see \cite{Ma, 
Theorem 1}).

\proclaim{Corollary 1.1} Let $0\ls a\ls k-1$, then
$$
f_{1,k-1}^a(n)=\left\{ \alignedat 2
(&k-2)!(k-1)^{n+2-k} &\quad &\text{for $n\gs k$}\\
&n!                  &\quad & \text {for $n<k$}.
\endalignedat\right.
$$
\endproclaim

\demo{Proof} Since $R_{1,k-1}=1+(k-1)x$, the Main Theorem implies
$$
F_{1,k-1}^a(x)=\frac{1-\sum_{r=0}^{k-3}x^r(k-r-2)r!}{1-(k-1)x}=
\frac{x^{k-2}(k-2)!}{1-(k-1)x}+\sum_{r=0}^{k-3}x^rr!,
$$
and the result follows.
\qed
\enddemo

Another immediate corollary of the Main Theorem gives the asymptotics for
$\flma(n)$ as $n\to\infty$.

\proclaim{Corollary 1.2} $\flma(n)\sim c\gamma^n$, where $c$ is a constant 
depending on $l$ and $m$, and $\gamma$ is the maximal root of 
$L_{\lambda}^{\mu-\lambda}${\rm;} in particular, 
$\gamma\ls k-2+\sqrt{1+4(l-1)(m-1)}$.
\endproclaim

\demo{Proof} Follows from standard results in the theory of rational 
generating functions (see e.g. \cite{St, vol.~1, Ch.~4}) and the fact that 
all the 
roots of Laguerre polynomials are simple (see \cite{Sz, Ch.~3.3}).
The upper bound on $\gamma$ is obtained in \cite{IL}.
\qed
\enddemo

\heading 2. Proofs
\endheading

First of all, we make the following simple, though useful observation.

\proclaim{Lemma 2.1} For any natural $a$, $l$, $m$, $n$ such that 
$1\ls a\ls l+m-1$ one has $\flma(n)=f_{m,l}^{m+l-a}(n)$.
\endproclaim

\demo{Proof} Denote by $\rho_n$ and $\varkappa_n$ the involutions $S_n\to S_n$
that take $(i_1,i_2,\dots,i_k)$ to $(i_k,\dots,i_2,i_1)$ ({\it reversal\/})
to $(n+1-i_1,n+1-i_2,\dots,n+1-i_n)$ ({\it complement\/}), respectively. 
It is easy to
see that for any $T\subset S_k$, the involutions $\rho_n$ and $\varkappa_n$
provide natural bijections between the sets $S_n(T)$ and $S_n(\rho_kT)$,
and between $S_n(T)$ and $S_n(\varkappa_kT)$, respectively. It remains to
note that $\rho_k\varkappa_k\sigma^a\plm=\sigma^{l+m-a}\plm$.
\qed
\enddemo

From now on we assume that $a\gs0$, $l\gs1$, $m\gs1$ are fixed, and denote
$b=n-m+a$. It follows from Lemma~2.1 that we may assume that $a\ls m$, 
and hence $b\ls n$. This means, in other words, that
$\tau\in S_k$ belongs to $\sigma^a\plm$ if and only if $(\tau_1,\dots,\tau_l)$
is a permutation of the numbers $a+1,\dots,a+l$. In what follows we usually
omit the indices $a$, $l$, $m$ whenever appropriate; for example, instead of
$\flma(n)$ we write just $f(n)$.

For any $n\gs k$ and any $d$ such that $1\ls d\ls n$, we denote by 
$g_n(i_1,\dots,i_d)=g_{n;l,m}^a(i_1,\dots,i_d)$ the number of permutations
$\pi\in S_n(\sigma^a\plm)$ such that $\pi_j=i_j$ for $j=1,\dots,d$. It is
natural to extend $g_n$ to the case $d=0$ by setting $g_n(\varnothing)=f(n)$.

The following properties of the numbers $g_n(i_1,\dots,i_d)$ can be deduced
easily from the definitions.

\proclaim{Lemma 2.2} {\rm (i)} Let $n\gs k$ and $1\ls i\ls n$, then
$$
g_n(\dots,i,\dots,i,\dots)=0.
$$

{\rm (ii)} Let $n \gs k$ and $a+1\ls i_j\ls b$ for $j=1,\dots,l$, then
$$
g_n(i_1,\dots,i_l)=0.
$$

{\rm (iii)} Let $n\gs k$, $1\ls r\ls d\ls l$, and $a+1\ls i_j\ls b$ for 
$j=1,\dots,d$, $j\ne r$, then
$$
g_n(i_1,\dots,i_d)= \left\{\alignedat2
 &g_{n-1}(i_1-1,\dots, i_{r-1}-1,i_{r+1}-1,
\dots, i_d-1)&\quad &\text{if $1\ls i_r\ls a$}\\
&g_{n-1}(i_1,\dots, i_{r-1},i_{r+1},\dots, i_d)&\quad &\text{if 
$b+1\ls i_r\ls n$}.
\endalignedat\right.
$$
\endproclaim

\demo{Proof} Property (i) is evident. Let us prove (ii). By (i), we may assume
that the numbers $i_1,\dots,i_l$ are distinct. Take an arbitrary $\pi\in S_n$
such that $\pi_j=i_j$ for $j=1,\dots,l$. Evidently, for any $r\ls a$ there
exists a position $j_r>l$ such that $\pi_{j_r}=r$; the same is true for any
$r\gs b+1$. Therefore, the restriction of $\pi$ to the positions
$1,2,\dots,l,j_1,j_2,\dots,j_a,j_{b+1},j_{b+2},\dots,j_n$ (in the
proper order) gives an occurrence of $\tau\in \sigma^a\plm$ in $\pi$. Hence,
$\pi\notin S_n(\sigma^a\plm)$, which means that $g_n(i_1,\dots,i_l)=0$.

To prove (iii), assume first that $1\ls r\ls a$. Let $\pi\in S_n$
and $\pi_j=i_j$ for $j=1,\dots,d$. We define $\pi^*\in S_{n-1}$ by
$$
\pi^*_j= \left\{\alignedat2
 &\pi_j-1&\quad &\text{for $1\ls j\ls r-1$},\\
 &\pi_{j+1}-1 &\quad &\text{for $j\gs r$ and $\pi_{j+1}>i_r$},\\
 &\pi_{j+1} &\quad &\text{for $j\gs r$ and $\pi_{j+1}<i_r$}.
\endalignedat\right.\tag1
$$

We claim that $\pi\in S_n(\sigma^a\plm)$ if and only if $\pi^*\in 
S_{n-1}(\sigma^a\plm)$. Indeed, the only if part is trivial, since any
occurrence of $\tau\in \sigma^a\plm$ in $\pi^*$ immediately gives rise
to an occurrence of $\tau$ in $\pi$. Conversely, any occurrence of $\tau$
in $\pi$ that does not include $i_r$  gives rise
to an occurrence of $\tau$ in $\pi^*$. Assume that there exists an occurrence 
of $\tau$ in $\pi$ that includes $i_r$. Since $r\ls d\ls l$, this occurrence
of $\tau$ contains $a$ entries that are situated to the right of $i_r$ and
are strictly less than $i_r$. However, the whole $\pi$ contains only $a-1$
such entries, a contradiction. It now follows from (1) that property (iii)
holds for $1\ls i_r\ls a$. The case $b+1\ls i_r\ls n$ is treated similarly
with the help of transformation $(\pi\in S_n)\mapsto (\pi^\circ\in S_{n-1})$
given by
$$
\pi^\circ_j= \left\{\alignedat2 
&\pi_j&\quad &\text{for $1\ls j\ls r-1$},\\
&\pi_{j+1}-1 &\quad &\text{for $j\gs r$ and $\pi_{j+1}>i_r$},\\
&\pi_{j+1} &\quad &\text{for $j\gs r$ and $\pi_{j+1}<i_r$}.
\endalignedat\right.
$$
\qed
\enddemo

Now we introduce the quantity that plays the crucial role in the proof of
the Main Theorem. For $n\gs k$ and $1\ls d\ls l$ we put
$$
A(n,d)=A_{l,m}^a(n,d)=\sum_{i_1,\dots,i_d=a+1}^{b}g_n(i_1,\dots,i_d).
$$
As before, this definition is extended to the case $d=0$ by setting
$$
A(n,0)=g_n(\varnothing)=f(n).
$$

\proclaim{Theorem 2.3} Let $n\gs k+1$ and $1\ls d\ls l-1$, then
$$
A(n,d+1)=a(n,d)-(m-d)A(n-1,d)-dA(n-1,d-1). \tag2
$$
\endproclaim

\demo{Proof} First of all, we introduce two auxiliary sums:
$$\align
B(n,d)=B_{l,m}^a(n,d)&=\sum_{i_1,\dots,i_d=a+1}^{b+1}g_n(i_1,\dots,i_d),\\
C(n,d)=C_{l,m}^a(n,d)&=\sum_{i_1,\dots,i_d=a}^{b}g_n(i_1,\dots,i_d);
\endalign
$$
once again, $B(n,0)=C(n,0)=f(n)$.

Let us prove three simple identities relating together the sequences
$\{A(n,d)\}$, $\{B(n,d)\}$, $\{C(n,d)\}$.

\proclaim{Lemma 2.4} Let $n\gs k$ and $1\ls d\ls l$, then\/{\rm:}
$$\gather
A(n,d)=A(n,d-1)-(m-a)B(n-1,d-1)-aC(n-1,d-1),\\
(m-a)A(n,d)=(m-a)B(n,d)-(m-a)dB(n-1,d-1),\\
aA(n,d)=aC(n,d)-adC(n-1,d-1).
\endgather
$$
\endproclaim

\demo{Proof} To prove the first identity, observe that by 
definitions and Lemma~2.2(iii) for the case $r=d$, one has
$$\multline
A(n,d-1)-A(n,d)=\sum_{i_1,\dots,i_{d-1}=a+1}^{b}
\sum_{i_d=1}^ng_n(i_1,\dots,i_d)-A(n,d)\\
=\sum_{i_1,\dots,i_{d-1}=a+1}^{b}\left(
\sum_{i_d=1}^ag_n(i_1,\dots,i_d)+
\sum_{i_d=b+1}^ng_n(i_1,\dots,i_d)\right)\\
=\sum_{i_1,\dots,i_{d-1}=a+1}^{b}\big(
ag_{n-1}(i_1-1,\dots,i_{d-1}-1)+(m-a)g_{n-1}(i_1,\dots,i_{d-1})\big)\\
=aC(n-1,d-1)+(m-a)B(n-1,d-1),
\endmultline
$$
and the result follows.

The second identity is trivial for $a=m$, so assume that $0\ls a\ls m-1$
and observe that by definitions and Lemma~2.2(ii) and (iii), one has
$$\multline
B(n,d)= \sum_{i_1,\dots,i_d=a+1}^{b}g_n(i_1,\dots,i_d)\\+
\sum_{j=1}^d\sum_{i_1,\dots,\hat\imath_j,\dots,i_d=a+1}^{b}
g_n(i_1,\dots,i_{j-1},b+1,i_{j+1}\dots,i_d)\\
= A(n,d)+\sum_{j=1}^d\sum_{i_1,\dots,\hat\imath_j,\dots,i_d=a+1}^{b}
g_{n-1}(i_1,\dots,i_{j-1},i_{j+1},\dots,i_d)\\
= A(n,d)+dB(n-1,d-1),
\endmultline
$$
and the result follows.

Finally, the third identity is trivial for $a=0$, so assume that $1\ls a\ls m$
and observe that by definitions and Lemma~2.2(ii) and (iii), one has
$$\multline
C(n,d)= \sum_{i_1,\dots,i_d=a+1}^{b}g_n(i_1,\dots,i_d)\\+
\sum_{j=1}^d\sum_{i_1,\dots,\hat\imath_j,\dots,i_d=a+1}^{b}
g_n(i_1,\dots,i_{j-1},a,i_{j+1}\dots,i_d)
= A(n,d)\\+\sum_{j=1}^d\sum_{i_1,\dots,\hat\imath_j,\dots,i_d=a+1}^{b}
g_{n-1}(i_1-1,\dots,i_{j-1}-1,i_{j+1}-1,\dots,i_d-1)\\
= A(n,d)+dC(n-1,d-1),
\endmultline
$$
and the result follows.
\qed
\enddemo

Now we can complete the proof of Theorem~2.3. Indeed, using twice the first
identity of Lemma~2.4, one gets 
$$\align
A(n,d+1)&=A(n,d)-(m-a)B(n,d)-aC(n-1,d-1),\\
dA(n-1,d)&=dA(n-1,d-1)-d(m-a)B(n-2,d-1)-daC(n-2,d-1).
\endalign
$$
Next, the other two identities of Lemma~2.4 imply
$$
\multline
A(n,d+1)-dA(n-1,d)=A(n,d)-dA(n-1,d-1)\\
-\big((m-a)B(n-1,d)-(m-a)dB(n-2,d-1)\big)
-\big(aC(n-1,d)-adC(n-2,d-1)\big)\\
=A(n,d)-dA(n-1,d-1)-(m-a)A(n-1,d)-aA(n-1,d),
\endmultline
$$
and the result follows.
\qed
\enddemo

The next result relates the sequence $\{A(n,d)\}$ to the sequence $\{f(n)\}$.

\proclaim{Theorem 2.5} Let $n\gs k$ and $1\ls d\ls l$, then
$$
A(n,d)=\sum_{j=0}^d(-1)^jj!\binom mj \binom dj f(n-j).
$$
\endproclaim

\demo{Proof} Let $D(n,d)=D_{l,m}^a(n,d)$ denote the right hand side of the 
above identity. We claim that for $n\gs k+1$ and $1\ls d\ls l-1$, $D(n,d)$
satisfies the same relation (2) as $A(n,d)$ does. Indeed, 
$$\multline
D(n-1,d)=\sum_{j=0}^d(-1)^jj!\binom mj\binom dj f(n-1-j)\\
        =\sum_{j=1}^{d+1}(-1)^j(j-1)!\binom m{j-1}\binom d{j-1} f(n-j)\\
          +(m-d)(-1)^{d+1}d!\binom md f(n-d-1),
\endmultline
$$
and
$$\multline
D(n-1,d-1)=\sum_{j=0}^{d-1}(-1)^jj!\binom mj\binom {d-1}j f(n-1-j)\\
          =\sum_{j=1}^{d}(-1)^j(j-1)!\binom m{j-1}\binom {d-1}{j-1} f(n-j),
\endmultline
$$
and hence
$$
\multline
D(n,d)-(m-d)D(n-1,d)-dD(n-1,d-1)\\
=f(n)+(m-d)(-1)^{d+1}d!\binom md f(n-d-1)\\
+\sum_{j=1}^{d}(-1)^jj!
\left(\binom mj\binom dj+\frac{m-d}j\binom m{j-1}\binom d{j-1}+
\frac dj\binom m{j-1}\binom {d-1}{j-1}\right)f(n-j)\\
=f(n)+ \sum_{j=1}^{d}(-1)^jj!\binom mj\binom {d+1}j f(n-j)+
(-1)^{d+1}(d+1)!\binom m{d+1}f(n-d-1)\\=D(n,d+1).
\endmultline
$$

It follows that $D(n,d)$ (as well as $A(n,d)$) are defined uniquely for
$n\gs k$ and $1\ls l \ls d$ by initial values $D(k,d)$, $D(n,0)$, and $D(n,1)$
($A(k,d)$, $A(n,0)$, and $A(n,1)$, respectively). It is easy to see that for
$n\gs k$ one has $A(n,0)=D(n,0)=f(n)$. Next, the first identity of Lemma~2.4
for $d=1$ gives
$$
A(n,1)=A(n,0)-(m-a)B(n-1,0)-aC(n-1,0)=f(n)-mf(n-1)\quad\text{for $n\gs k$}.
$$
On the other hand, by definition,
$$
D(n,1)=f(n)-mf(n-1) 
\quad\text{for $n\gs k$},
$$
and hence $A(n,1)=D(n,1)$. Finally, a simple combinatorial argument shows that
$$
A(k,d)=d!\binom ld (k-d)! -l!m! \quad\text{for $1\ls d\ls l$}.
$$
On the other hand,
$$
D(k,d)=\sum_{j=0}^{d}(-1)^jj!\binom mj\binom dj (k-j)! -l!m!,
$$
since $f(r)=r!$ for $1\ls r\ls k-1$ and $f(k)=k!-l!m!$. To prove $A(k,d)=
D(k,d)$ it remains to check that
$$
\sum_{j=0}^{d}(-1)^jj!\binom mj\binom dj (k-j)!=d!\binom ld (k-d)!,
$$
which follows from Lemma~2.6 below.
\qed
\enddemo

Finally, we are ready to prove the Main Theorem stated in Sec.~1. 
First of all, by Lemma~2.2(ii), $A(n,l)=0$ for $n\gs k$. Hence, by Theorem~2.5,
$$
\sum_{j=0}^l(-1)^jj!\binom mj\binom lj f(n-j)=0 \quad\text{for $n\gs k$},
$$
or, equivalently,
$$
\sum_{j=0}^l(-1)^jj!\binom mj\binom lj x^j f(n-j)x^{n-j}=0 
\quad\text{for $n\gs k$}.
$$
As it was already mentioned, $f(r)=r!$ for $1\ls r\ls k-1$, therefore
$$
\sum_{j=0}^l(-1)^jj!\binom mj\binom lj x^j \left(\Flma(x)
-\sum_{i=0}^{k-j-1}x^ii!\right)=0.\tag3
$$

Recall that the rook polynomial of the rectangular $s\times t$ board, 
$s\ls t$, satisfies relation
$$
R_{s,t}(x)=\sum_{j=0}^sj!\binom tj \binom sj x^j
$$
(see \cite{Ri, Ch.~7.4}). Hence, (3) is equivalent to
$$\gather
\Flma(x)R_{\lambda,\mu}(-x)=\sum_{j=0}^l(-1)^jj!\binom mj\binom lj x^j
\sum_{i=0}^{k-j-1}x^ii!\\
=\sum_{r=0}^{k-1}x^rr!\sum_{j=0}^r(-1)^j\frac{\binom lj\binom mj}{\binom rj}\\
=\sum_{r=0}^{\lambda-1}x^rr!
\sum_{j=0}^r(-1)^j\frac{\binom lj\binom mj}{\binom rj}
+\sum_{r=\lambda}^{\mu-1}x^rr!
\sum_{j=0}^\lambda(-1)^j\frac{\binom lj\binom mj}{\binom rj}+
\sum_{r=\mu}^{k-1}x^rr!
\sum_{j=0}^\lambda(-1)^j\frac{\binom lj\binom mj}{\binom rj}.
\endgather
$$
By Lemma 2.6 below, the third term of the above expression vanishes, while
the second term is equal to
$$
\sum_{r=\lambda}^{\mu-1}x^rr!(-1)^\lambda\frac{(r-\lambda)!(k-r-1)!}
{\mu-r-1)!r!},
$$
and the first expression of the Main Theorem follows. The second expression
is obtained easily from (3) and relation between rook polynomials and Laguerre
polynomials given in Sec.~1.
\qed

It remains to prove the following technical result, which is apparently known;
however, we failed to find a reference to its proof, and decided to present
a short proof inspired by the brilliant book \cite{PWZ}.

\proclaim{Lemma~2.6} Let $1\ls s\ls t$ and let 
$$
M(s,t)=\sum_{i=0}^s(-1)^j\frac{\binom si\binom ti}{\binom ni}.
$$
Then\/{\rm :}
$$
M(s,t)=\left\{\alignedat 2
 &\frac{\binom {n-t}s}{\binom ns}\quad &\text{if $n\gs s+t$},\\
              &0\quad &\text{if $t\ls n\ls s+t-1$},\\
   &(-1)^s\frac{\binom {s+t-n-1}s}{\binom ns}\quad &\text{if $s\ls n\ls t-1$}.
\endalignedat \right.
$$
\endproclaim

\demo{Proof} Direct check reveals that $M(s,t)$ is a hypergeometric series;
to be more precise,
$$
M(s,t)={}_2F_1\left[^{-t,-s}_{-n}; 1\right].
$$
Since $-s$ is a nonpositive integer, the Gauss formula applies 
(see \cite{PWZ, Ch.~3.5}), and we get
$$
M(s,t)=\lim_{z\to n}\frac{\Gamma(-z+t+s)\Gamma(-z)}{\Gamma(t-z)\Gamma(s-z)}.
$$
Recall that 
$$
\Gamma(x)\Gamma(1-x)=\frac\pi{\sin\pi x}.\tag4
$$

If $n\gs s+t$, we apply (4) for $x=-z+t+s$, $x=t-z$, $x=s-z$, $x=-z$, and
get
$$
M(s,t)=-\frac{\Gamma(n-t+1)\Gamma(n-s+1)}{\Gamma(n-t-s+1)\Gamma(n+1)}
\lim_{z\to n}\frac{\sin\pi(t-z)\sin\pi(s-z)}{\sin\pi z\sin\pi(t+s-z)}=
\frac{\binom {n-t}s}{\binom ns}.
$$

If $t\ls n\ls s+t-1$, we apply (4) for  $x=t-z$, $x=s-z$, $x=-z$, and
get
$$
M(s,t)=-\frac{\Gamma(n-t+1)\Gamma(n-s+1)\Gamma(s+t-n)}{\Gamma(n+1)}
\lim_{z\to n}\frac{\sin\pi(t-z)\sin\pi(s-z)}{\sin\pi z}=0.
$$

Finally, if $s\ls n\ls t-1$, we apply (4) for  $x=s-z$, $x=-z$, and
get
$$
M(s,t)=-\frac{\Gamma(s+t-n)\Gamma(n-s+1)}{\Gamma(t-n)\Gamma(n+1)}
\lim_{z\to n}\frac{\sin\pi(s-z)}{\sin\pi z}=
(-1)^s\frac{\binom {s+t-n-1}s}{\binom ns}.
$$
\qed
\enddemo

\heading 3. Concluding remarks
\endheading

Observe first, that according to the Main Theorem, $\Flma(x)$ does not depend
on $a$; in other words, $|S_n(\plm)|=|S_n(\sigma^a\plm)|$ for any $a$. We
obtained this fact as a consequence of lengthy computations. A natural
question would be to find a bijection between $S_n(\plm)$ and 
$S_n(\sigma^a\plm)$ that explains this phenomenon.

Second, it is well known that rook polynomials (or the corresponding Laguerre
polynomials) are related to permutations with restricted positions, see
\cite{Ri, Ch.7,8}. Laguerre polynomials also arise in a natural way in the 
study of generalized derangements (see \cite{FZ} and references therein).
It is tempting to find a combinatorial relation between permutations with 
restricted positions and permutations avoiding maximal parabolic subgroups,
which could explain the occurrence of Laguerre polynomials in the latter
context.

Finally, one can consider permutations avoiding nonmaximal parabolic subgroups
of $S_k$. The first natural step would be to treat the case of subgroups
generated by $k-3$ simple transpositions. It is convenient
to denote by $P_{l_1,l_2,l_3}$ (with $l_1+l_2+l_3=k$) the subgroup of $S_k$
generated by all the simple transpositions except for $s_{l_1}$ and
$s_{l_1+l_2}$; further on, we set $f_{l_1,l_2,l_3}(n)=|S_n(P_{l_1,l_2,l_3})|$,
and $F_{l_1,l_2,l_3}(x)=\sum_{n\gs0}f_{l_1,l_2,l_3}(n)x^n$. It is easy to see
that $F_{l_1,l_2,l_3}(x)=F_{l_3,l_2,l_1}(x)$, so one can assume that
$l_1\ls l_3$. This said,
the main result of \cite{BDPP} can be formulated as follows: let $k\gs3$,
then
$$
F_{1,1,k-2}(x)=\sum_{r=1}^{k-2}x^rr!+\frac{(k-3)!x^{k-4}}2
\left(1-(k-1)x-\sqrt{1-2(k-1)x+(k-3)^2x^2}\right).
$$
To the best of our knowledge, this is the only known instance of
$F_{l_1,l_2,l_3}(x)$. It is worth to note that even in this, simplest, case 
of nonmaximal parabolic subgroup, the generating function is no more rational.

\enddocument